\documentclass[12pt]{article}
\newtheorem{theorem}{Theorem}[section]
\newtheorem{proposition}[theorem]{Proposition}
\newtheorem{lemma}[theorem]{Lemma}
\newtheorem{corollary}[theorem]{Corollary}

\newcommand{\nl}{\newline}

\newcommand{\R}{{\bf R}}

\newcommand{\cD}{{\cal D}}

\newcommand{\cW}{{\cal W}}
\newcommand{\cH}{{\cal H}}

\newcommand{\diver}{{\rm div}}

\newcommand{\parder}[2]{\frac{\partial{#1}}{\partial{#2}}}

\newcommand{\darr}[4]{{\left\{\begin{array}{ll}
   {#1}&{#2}\\
   {#3}&{#4}
\end{array}\right.}}

\newcommand{\darrn}[4]{{\begin{array}{ll}
   {#1}&{#2}\\
   {#3}&{#4}
\end{array}}}

 \newcommand{\ia}{({\rm i})}
 \newcommand{\ib}{({\rm ii})}
 \newcommand{\ic}{({\rm iii})}

 \title{ Critical heat kernel estimates via    Hardy-Sobolev inequalities}
 \author{
G. Barbatis\footnote{Department of Mathematics,
 University of Ioannina, 45110 Ioannina, Greece}
 \and S. Filippas\footnote{Department of Applied Mathematics,
 University of Crete, 71409 Heraklion, Greece} \and A. Tertikas
 \footnote{Department of Mathematics,
 University of Crete, 71409 Heraklion, Greece and  \nl
Institute of  Applied and Computational Mathematics,
FORTH, 71110 Heraklion, Greece}
}
\date{\today}

\begin{document}
\maketitle


\begin{abstract}
We obtain Sobolev inequalities for the Shcr\"{o}dinger operator $-\Delta-V$, where
$V$ has critical behaviour $V(x)=((N-2)/2)^2|x|^{-2}$ near the origin. We apply
these inequalities to obtain pointwise estimates on the associated heat kernel,
improving upon earlier results.

\end{abstract}

\section{Introduction}

The purpose of this paper is to obtain some new Hardy-Sobolev inequalities
and then use them in order to obtain new  heat kernel estimates for the Schr\"{o}dinger operator
$-\Delta -V$ for positive potentials $V$ with critical singularities, improving
upon analogous estimates of this type.

As a typical example, let us consider the case of a bounded domain
$\Omega\subset\R^N$, $N\geq 3$, containing the origin.
We obtain upper estimates on the heat kernel of the operator
\begin{equation}
Hu =-\Delta u-\frac{\lambda}{|x|^2}u ,\quad u|_{\partial\Omega}=0 ,
\label{eq:i1}
\end{equation}
for various values of the real parameter $\lambda$. It is well known
that the power $|x|^{-2}$ is critical and that the associated
heat kernel exhibits behaviour which is different from that
of the case $|x|^{-\gamma}$, $\gamma<2$, which is in the  Kato class.
Recently Milman and Semenov \cite{MS} obtained upper heat kernel estimates
for the operator (\ref{eq:i1}) when $\lambda<((N-2)/2)^2$.
In the case $0<\lambda<((N-2)/2)^2$ they showed that,
modulo the usual Gaussian term, the heat kernel satisfies
\[ K(t,x,y)<ct^{-\frac{N}{2}}|x|^{-\alpha}|y|^{-\alpha},\]
where $\alpha$ denotes the smallest solution of $\alpha(N-2-\alpha)=\lambda$.
In Theorem \ref{thm:mainthm} we extend this result to the critical case
$\lambda=((N-2)/2)^2$, namely we prove that the corresponding kernel satisfies
\[ K(t,x,y)<ct^{-N/2}|x|^{-\frac{N-2}{2}}|y|^{-\frac{N-2}{2}}.\]
This estimate is sharp as can be seen by comparing with the results
of V\'{a}zquez and Zuazua \cite{VZ}.
Moreover, for the subcritical case $\lambda<((N-2)/2)^2$ and for $\lambda<0$
we improve the small-time dependence obtained in \cite{MS}, namely we prove the
estimate
$K(t,x,y)<ct^{-\frac{N}{2}+ \alpha}|x|^{-\alpha}|y|^{-\alpha}$ instead of
$K(t,x,y)<ct^{-\frac{N}{2}+\alpha   -\epsilon}|x|^{-\alpha}|y|^{-\alpha}$; see
Proposition \ref{prop:pos}.

In addition we  consider operators that act on the whole of $\R^N$ with potentials
having the critical Hardy singularity near zero, of the form
\begin{equation}
V_{\epsilon}(x)=\darr{\Bigl(\frac{N-2}{2}\Bigr)^2|x|^{-2},}{|x|<1,}{\epsilon f(x),}{|x|>1,}
\label{eq:tax}
\end{equation}
under appropriate subcritical  assumptions on the positive function $f$. Thus, in  Theorem \ref{thm:elec}
it is shown that if $\epsilon>0$ is small enough then
the heat kernel of $-\Delta-V_{\epsilon}$ satisfies
\begin{equation}
K(t,x,y)<ct^{-N/2}\max\{ |x|^{-\frac{N-2}{2}},1\}\max\{ |y|^{-\frac{N-2}{2}},1\}.
\label{eq:fire}
\end{equation}

We also consider potentials that exhibit the critical behaviour $((N-2)/2)^2|x|^{-2}$
near infinity, that is
\begin{equation}
\hat{V}_{\epsilon}(x)=\darr{\epsilon g(x),}{|x|<1,}{(\frac{N-2}{2})^2|x|^{-2},}{|x|>1.}
\label{eq:aek}
\end{equation}
Under appropriate  subcritical assumptions on $g$ we obtain Sobolev estimates for
$-\Delta-\hat{V}_{\epsilon}$ for a sharp range of $\epsilon>0$.
We note here that while the question of Sobolev inequalities for $-\Delta-\hat{V}_{\epsilon}$
is rather similar to that for $-\Delta-V_{\epsilon}$, when it comes
to heat kernel estimates essential differences arise. As mentioned earlier,
the Sobolev inequality
for $-\Delta-V_{\epsilon}$ yields estimate (\ref{eq:fire}) for the corresponding heat kernel.
On the other hand, while the short-time behaviour of the heat kernel of
$-\Delta-\hat{V}_{\epsilon}$ is similar to that of the Laplacian, the long-time behaviour
is very different. Zhang \cite{Z} used a parabolic Harnack inequality to
obtain precise estimates for the heat kernel of
$-\Delta-V$, when $V$ is equal to $\lambda|x|^{-2}$, $\lambda<((N-2)/2)^2$,
near infinity.
This generalizes the   estimates given by Davies and Simon \cite{DS}, where, however, more precise
estimates were given. The corresponding problem for the critical
case $\lambda=((N-2)/2)^2$ remains open.

Going back to bounded $\Omega\subset\R^N$ and to the operator $H$ defined in (\ref{eq:i1}),
for the critical case $\lambda=((N-2)/2)^2$ we finally  consider additional
 singularities, that is, we consider potentials of the form
$((N-2)/2)^2|x|^{-2}+V_1$, where $V_1>0$ is also critical; $V_1$ is defined
as a series involving iterated logarithms (see definition (\ref{eq:anti})) and is critical
 in the sense that the following improved Hardy inequality holds
\begin{equation}
 \int_{\Omega}|\nabla u|^2 - \Bigl(\frac{N-2}{2}\Bigr)^2\int_{\Omega}\frac{u^2}{|x|^2} \geq
\int_{\Omega}V_1u^2,\quad u\in C^{\infty}_c(\Omega),
\label{1.5}
\end{equation}
whereas this inequality is no longer true if we replace $V_1$ by $(1+\epsilon)V_1$ for any $\epsilon>0$.
It is remarkable that the extra potential $V_1$ does not affect the time dependence of the heat kernel
 estimates,
but only affects the spatial singularity at the origin (cf Theorems  \ref{thm:mainthm}  and  \ref{thm:lognon}).
This is in contrast with Proposition  \ref{prop:pos}(ii)  where,
 for $\lambda<0$, the  potential affects the 
 the time singularity of the heat kernel  as well.

Throughout the paper we study a number of concrete potentials. These are chosen
precisely because they are critical. By simple monotonicity one can then obtain heat kernel estimates for
a whole range of other potentials, including potentials that are not radially symmetric.

To prove the above heat kernel estimates we first use 
 an appropriate change of variables, $u=\phi w$, by means of which, 
 the problem is reduced to obtaining
uniform estimates on the heat kernel $K_{\phi}(t,x,y)$ of an auxiliary operator $H_{\phi}$
which acts on the function $w$; see, e.g., [MS]. Those estimates are in turn proved by means of some
new Hardy- Sobolev inequalities.

As a typical example of such an inequality we mention the following inequality proved 
by Brezis and V\'{a}zquez [BV]:
\begin{equation}
 \int_{\Omega}|\nabla u|^2 - \Bigl(\frac{N-2}{2}\Bigr)^2\int_{\Omega}\frac{u^2}{|x|^2} \geq
K \Bigl(\int_{\Omega}|u|^p \, dx \Bigr)^{2/p},
\label{1.6}
\end{equation}
valid for $u \in H_0^1(\Omega)$ and $1<p<\frac{2N}{N-2}$; this inequality fails for the critical
Sobolev exponent $p= \frac{2N}{N-2}$.  To obtain sharp heat kernel estimates one needs to go up
to the critical exponent. In connection with this we mention the following sharp Hardy-Sobolev
inequality established in [FT]:
\begin{equation}
 \int_{\Omega}|\nabla u|^2 - \Bigl(\frac{N-2}{2}\Bigr)^2\int_{\Omega}\frac{u^2}{|x|^2} \geq
c \Bigl( \int_{\Omega} |u|^{\frac{2N}{N-2}} X_1^{1+\frac{N}{N-2}}
 \Bigl( \frac{|x|}{D} \Bigr)dx \Bigr)^{\frac{N-2}{N}},
\label{1.7}
\end{equation}
valid for $u \in H_0^1(\Omega)$; here $D=\sup_{\Omega}|x|$ and $X_1(t)=(1-\log t)^{-1}$, $t \in (0,1)$.
In the present work we derive new  Hardy-Sobolev inequalities that involve potentials such as the ones
given in  (\ref{eq:tax})   or   (\ref{eq:aek}); see Theorems \ref{thm:gian}, \ref{thm:newsob}, 
 \ref{thm:gian1} and \ref{thm:newsob1}.  We should  mention  that the validity of improved Hardy
inequalities is strongly connected to the existence and large time behaviour of solutions of the
heat equation with singular potential; see, e.g.,  Brezis and V\'{a}zquez [BV], Cabr\'{e} and Martel  [CM],
D\'{a}vila and Dupaigne [DD]
as well as   V\'{a}zquez and Zuazua [VZ].

As a byproduct of our approach we establish various results concerning improved Hardy inequalities
with boundary terms.  Such inequalities have recently attracted attention, see  Adimurthi and Esteban [AE],
Wang and Zhu [WZ] and references therein.

The structure of the paper is as follows: in Section 2 we present some
auxiliary results  concerning improved Hardy inequalities
with boundary terms.  In Section 3
we prove the Hardy-Sobolev inequalities; in Section 4 we apply them to obtain
heat kernel estimates. Finally, in Section 5 we prove refined Sobolev inequalities
and heat kernel estimates when additional singularities are present.

{\bf Acknowledgment} We acknowledge partial  support by the  RTN  European network
 Fronts--Singularities, HPRN-CT-2002-00274.

\setcounter{equation}{0}
\section{Two minimization problems}
Throughout this section, $\Omega  \subset  \R^N$, $N\geq 3$,  is a bounded domain containing
the origin with $C^1$ boundary. Also, we always  denote   by $\nu$ the outward-pointing
(with respect to $\Omega$) unit vector on the surface $\partial \Omega$.
In Section 2.1 we will work on $\Omega$, while in Section 2.2 we will work on $\Omega^c$.
The results of this section will be applied in Section 3.

\subsection{Bounded domains}

For $ \alpha>0$ we define
\begin{equation}
\lambda_{\Omega}(\alpha) = 
\inf_{H^1(\Omega)}
{\frac{\int_{\Omega}|\nabla u|^2dx +  \alpha\int_{\partial \Omega}   \frac{x \cdot \nu}{|x|^2}
  u^2 dS}{\int_{\Omega} \frac{u^2}{|x|^2} dx}}.
\label{2.2.1}
\end{equation}

\begin{lemma}
Let $\Omega$ be a bounded domain in $\R^N$, $N\geq 3$, containing the origin. \nl
$\ia$ If  $0< \alpha\leq \frac{N-2}{2}$, then $\lambda_{\Omega}(\alpha)
 = \alpha(N-2-\alpha)$. Moreover,
$|x|^{-\alpha} \in H^1(\Omega)$ is a minimizer for $0<\alpha  <\frac{N-2}{2}$, whereas  for
$\alpha =\frac{N-2}{2}$ there is no $H^1(\Omega)$ minimizer. \nl
$\ib$ If $\alpha >\frac{N-2}{2}$ and $\Omega$ is starshaped with respect to zero, then 
$\lambda_{\Omega}(\alpha ) = \Bigl(\frac{N-2}{2} \Bigr)^2$ and there is
no $H^1(\Omega)$ minimizer.
\label{lem:intf}
\end{lemma}

In case  $\Omega$ is not starshaped with respect to zero, concerning the analogue of
part (ii) of the above Lemma,
we have 
\begin{lemma}
Let $\Omega$ be a bounded domain in $\R^N$, $N\geq 3$, containing the origin,
which is not starshaped with respect to zero. Then, there exist
 finite constants $\alpha^* \geq N-2$
 and $\alpha_*  \in [\frac{N-2}{2}, \alpha^*)$ depending on $\Omega$ such that: \nl
$\ia$ $\lambda_{\Omega}(\alpha^*)=0$, whereas $\lambda_{\Omega}(\alpha )>0$ for all
$\frac{N-2}{2} < \alpha <\alpha^*$. \nl
$\ib$ If $\frac{N-2}{2} \leq  \alpha  \leq \alpha_*$, then 
 $\lambda_{\Omega}(\alpha ) = \bigl(\frac{N-2}{2} \big)^2$ and 
 there is no $H^1(\Omega)$ minimizer. \nl
$\ic$ If  $\alpha_*  < \alpha < \alpha^*$, then  $\max(0,\,  \alpha(N-2-\alpha))
 \leq \lambda_{\Omega}(\alpha)
 < \big(\frac{N-2}{2} \big)^2$,  and there exists an $H^1(\Omega)$ minimizer.
\label{lem:intf2}
\end{lemma}
{\bf Remarks. 1.} We note  in particular that for any $\Omega$ and any $\alpha>0$ there holds
\begin{equation}
\alpha(N-2-\alpha) \leq  \lambda_{\Omega}(\alpha) \leq  \Bigl(\frac{N-2}{2} \Bigr)^2.
\label{eq:rem}
\end{equation}
{\bf 2.} We do not know whether there exists a non-starshaped domain $\Omega$ with
smooth boundary so that $\alpha^*=N-2$. Similarly, we do not know whether there there
exists such an $\Omega$ for which $\alpha_*=(N-2)/2$.

{\em Proof of Lemmas \ref{lem:intf} and \ref{lem:intf2}:}
 Let $u \in C^{\infty}(\overline{\Omega})$ be supported outside a neighborhood of zero. 
 For any
$\alpha>0$ we set $u(x) = |x|^{-\alpha} v(x)$. A straightforward calculation shows that
\begin{equation}
\int_{\Omega}|\nabla u|^2dx = \int_{\Omega} |x|^{-2 \alpha} |\nabla v|^2dx
+  \alpha(N-2-\alpha) \int_{\Omega} \frac{u^2}{|x|^2}dx - \alpha
\int_{\partial \Omega} \frac{x \cdot \nu}{|x|^2}
  u^2 dS,
\end{equation}
therefore,
\begin{equation}
\int_{\Omega}|\nabla u|^2dx +\alpha
\int_{\partial \Omega} \frac{x \cdot \nu}{|x|^2}
  u^2 dS  \geq  \alpha(N-2-\alpha) \int_{\Omega} \frac{u^2}{|x|^2}dx.
\label{eq:22.1}
\end{equation}
By a simple density argument this inequality is valid for all $u \in  H^1(\Omega)$. This
implies in particular  the lower bound on $\lambda_{\Omega}(\alpha)$ in (\ref{eq:rem}).

If   $0< \alpha < \frac{N-2}{2}$, then $|x|^{-\alpha}$ is in  $H^1(\Omega)$ and an easy
calculation shows that satisfies (\ref{eq:22.1}) as equality, hence it is a minimizer.
If $\alpha =\frac{N-2}{2}$ then the fact that         
$ \lambda_{\Omega}(\alpha)= \Bigl(\frac{N-2}{2} \Bigr)^2$
follows by considering  the minimizing sequence $u_{\epsilon}(x) = |x|^{-\frac{N-2}{2}
+ \epsilon} \in H^1(\Omega)$, $\epsilon \rightarrow 0^+$. 

By the  same  sequence  $u_{\epsilon}(x) = |x|^{-\frac{N-2}{2}
+ \epsilon}  \in H^1(\Omega)$ one can show that
$ \Bigl(\frac{N-2}{2} \Bigr)^2 \geq   \lambda_{\Omega}(\alpha)$, 
for any $\alpha >0$ and any $\Omega$, thus proving the upper bound in  (\ref{eq:rem}).

Suppose now that $\Omega$ is starshaped and  $\alpha >\frac{N-2}{2}$. Then,
using first the fact that   $x \cdot \nu \geq 0$
on the boundary of $\Omega$,  and then  (\ref{eq:22.1}) (with $\alpha=\frac{N-2}{2}$)
\begin{eqnarray*}
\int_{\Omega}|\nabla u|^2dx +\alpha
\int_{\partial \Omega} \frac{x \cdot \nu}{|x|^2}   u^2 dS  &  \geq &
\int_{\Omega}|\nabla u|^2dx +\frac{N-2}{2}
\int_{\partial \Omega} \frac{x \cdot \nu}{|x|^2} u^2 dS \\
& \geq & \Bigl(\frac{N-2}{2}\Bigr)^2 \int_{\Omega} \frac{u^2}{|x|^2} dx.
\end{eqnarray*}
Hence, in this case  $\lambda_{\Omega}(\alpha )= \bigl(\frac{N-2}{2} \big)^2$.

We next show that when $\lambda_{\Omega}(\alpha ) =  \bigl(\frac{N-2}{2} \big)^2$,
there is no $H^1(\Omega)$ minimizer.  Indeed assuming that there is one, then 
it would be a positive   $H^1(\Omega)$ solution of the Euler-Lagrange equation
$$
\Delta u + \frac{\Bigl(\frac{N-2}{2}\Bigr)^2}{|x|^2}u =0, \qquad x \in \Omega,
$$
under  suitable Robin boundary conditions. However, this equation has no  $H^1(\Omega)$
positive solutions; see e.g. [FT, Theorem C] for a more general statement. (We note that
although in  this Theorem Dirichlet condition were imposed, the proof is independent
of the boundary conditions.) Thus, Lemma \ref{lem:intf} has been proved.

Suppose now that $\Omega$ is not starshaped with respect to zero.
The existence of  $\alpha^*$ 
follows from the continuity of  $\lambda_{\Omega}(\alpha)$ with respect to $\alpha$ 
combined with the fact that if $\Omega$ is not starshaped with respect to zero, 
then one can easily find test functions making the surface integral in (\ref{2.2.1})
negative. The fact that $\alpha^* \geq N-2$ follows from the lower bound in (\ref{eq:rem}).

From Lemma  \ref{lem:intf}(i), we have that 
 $\lambda_{\Omega}(\frac{N-2}{2}  ) = \bigl(\frac{N-2}{2} \big)^2$. We  then define 
$\alpha_*$  as the supremum of all $\alpha$ for which $\lambda_{\Omega}(\alpha)=
 \bigl(\frac{N-2}{2} \big)^2$. Assuming that $\alpha_*>\frac{N-2}{2}$, we will show that
for any $\frac{N-2}{2}< \alpha<\alpha_*$ there holds 
 $\lambda_{\Omega}(\alpha)=
 \bigl(\frac{N-2}{2} \big)^2$.
Indeed, if this is not the case then there would exist an $\alpha$ in the above interval
and $\phi \in  H^1(\Omega)$ such that
\begin{equation}
{\frac{\int_{\Omega}|\nabla \phi|^2dx +  \alpha\int_{\partial \Omega}   \frac{x \cdot \nu}{|x|^2}
  \phi^2 dS}{\int_{\Omega} \frac{\phi^2}{|x|^2} dx}} <\Big(\frac{N-2}{2} \Big)^2
\nonumber
\end{equation}
On the other hand from Lemma \ref{lem:intf}(i), we have that 
\begin{equation}
{\frac{\int_{\Omega}|\nabla \phi|^2dx + (\frac{N-2}{2})  \int_{\partial \Omega}
\frac{x \cdot \nu}{|x|^2}
  \phi^2 dS}{\int_{\Omega} \frac{\phi^2}{|x|^2} dx}} \geq \Big(\frac{N-2}{2} \Big)^2.
\nonumber
\end{equation}
From the above two inequalities it follows that $ \int_{\partial \Omega}
\frac{x \cdot \nu}{|x|^2}\phi^2 dS <0$.  Using $\phi$ as a test function and the fact that
$\alpha_*>\alpha$  we conclude that 
$\lambda_{\Omega}(\alpha_*)<(\frac{N-2}{2} \big)^2$, which is a contradiction.
Thus, the estimates of part (ii) and (iii) of Lemma \ref{lem:intf2} have been proved.

The nonexistence of $ H^1(\Omega)$ minimizer of part (ii) follows exactly as in Lemma 
\ref{lem:intf}.  The existence of $ H^1(\Omega)$ minimizer of part (iii)
will follow  later from a more general result; see Proposition \ref{lem:2.3}. $\hfill //$

\subsection{Complement of bounded domains}

Here we consider the complement of a bounded domain and we study the corresponding
infimum, that is
\begin{equation}
\mu_{\Omega}(\alpha) = \inf_{u\in  C^{\infty}_c(\R^N)|_{\Omega^c}}             
{\frac{\int_{\Omega^c}|\nabla u|^2dx -  \alpha\int_{\partial \Omega}   \frac{x \cdot \nu}{|x|^2}
  u^2 dS}{\int_{\Omega^c} \frac{u^2}{|x|^2} dx}}.
\label{2.2.2}
\end{equation}
where $\Omega$, as before, is a  bounded domain containing the origin
and  $\nu$ is  the outward-pointing (with respect to $\Omega$)
unit vector on the surface $\partial \Omega$.
Also,  $C^{\infty}_c(\R^N)|_{\Omega^c}$ is  the set of restrictions on $\Omega^c$ of
 all functions
$u\in C^{\infty}_c(\R^N)$. We also introduce the following norms,
\begin{eqnarray}
\|u\|_{\cD^{1,2}(\Omega^c)} & = &  (\int_{\Omega^c}|\nabla u|^2dx)^{1/2} +
 (\int_{\Omega^c}|u|^{\frac{2N}{N-2}} dx )^{\frac{N-2}{2N}} \\
\|u\|_{{\cal H}^1(\Omega^c)}  & = & (\int_{\Omega^c}|\nabla u|^2dx)^{1/2} +
 (\int_{\Omega^c}\frac{|u|^2}{|x|^2}dx )^{1/2}  \\
\|u\|_{{\cal W}(\Omega^c)}  & = & (\int_{\Omega^c}|\nabla u|^2dx)^{1/2} +
(\int_{\partial \Omega} |u|^{\frac{2(N-1)}{N-2}} dS)^{\frac{N-2}{2(N-1)}}
\end{eqnarray}
and we denote by $\cD^{1,2}(\Omega^c)$, ${\cal H}^{1}(\Omega^c)$ and ${\cal W}(\Omega^c)$
the completion of $C^{\infty}_c(\R^N)|_{\Omega^c}$ under the corresponding norms.
The space $\cW(\Omega^c)$ is well studied, see e.g. \cite{M}. For our
purposes however, the natural spaces to use are $\cD^{1,2}(\Omega^c)$ and $\cH^1(\Omega^c)$.
In the next lemma we show that these three spaces coincide (a trivial fact if $\Omega^c$
were replaced by $\Omega$).
\begin{lemma}
Let $\Omega$ be a bounded domain in $\R^N$, $N\geq 3$, with $C^1$ boundary,
containing the origin. Then $\cD^{1,2}(\Omega^c)={\cal H}^{1}(\Omega^c)={\cal W}(\Omega^c)$.
\end{lemma}
{\em Proof:} We will show that all norms are equivalent.
Let $u \in C^{\infty}_c(\R^N)|_{\Omega^c}$.  Under our assumptions, it follows easily
as in  Lemma \ref{lem:intf}  (cf (\ref{eq:22.1}) with  $\alpha = \frac{N-2}{2}$) that
\begin{equation}
\int_{\Omega^c}|\nabla u|^2dx - \frac{N-2}{2}
\int_{\partial \Omega} \frac{x \cdot \nu}{|x|^2}
  u^2 dS  \geq   \bigl(\frac{N-2}{2} \big)^2   \int_{\Omega^c} \frac{u^2}{|x|^2}dx,  \nonumber 
\end{equation}
whence,
\begin{eqnarray}
\int_{\Omega^c}\frac{|u|^2}{|x|^2}dx  & \leq & C \Bigl(\int_{\Omega^c}|\nabla u|^2dx
+ \int_{\partial \Omega} |u|^{2} dS  \Bigr)  \nonumber \\
& \leq & C  \Bigl(\int_{\Omega^c}|\nabla u|^2dx +
\Bigl(\int_{\partial \Omega} |u|^{\frac{2(N-1)}{N-2}} dS \Bigr)^{\frac{N-2}{(N-1)}}\Bigl).
\nonumber 
\end{eqnarray}
Hence, $\|u\|_{{\cal H}^{1}(\Omega^c)} \leq C \|u\|_{{\cal W}(\Omega^c)}$.
To obtain the  reverse inequality we note that
it follows  from the standard trace Theorem (e.g [A, Theorem 5.22], Chapter V)
-- applied to  $B \setminus \Omega$ for some  ball    $B \supset \Omega$ -- that
\begin{eqnarray}
 \Bigl(\int_{\partial \Omega}
 |u|^{\frac{2(N-1)}{N-2}} dS\Bigr)^{\frac{N-2}{2(N-1)}} \leq C  \|u\|_{{\cal H}^{1}(\Omega^c)}
 \nonumber   \\
 \Bigl(\int_{\partial \Omega}
 |u|^{\frac{2(N-1)}{N-2}} dS\Bigr)^{\frac{N-2}{2(N-1)}} \leq  C \|u\|_{\cD^{1,2}(\Omega^c)} 
 \nonumber  
\end{eqnarray}
From the first one it follows that $\|u\|_{{\cal W}(\Omega^c)}
 \leq C\|u\|_{{\cal H}^1(\Omega^c)} $, whence ${\cal H}^1(\Omega^c)={\cal W}(\Omega^c)$. From
the second one it follows that  $\|u\|_{{\cal W}(\Omega^c)}  \leq  C \|u\|_{\cD^{1,2}
(\Omega^c)}$.
Thus, it remains to prove that $ \|u\|_{\cD^{1,2}(\Omega^c)}  \leq  
C \|u\|_{{\cal W}(\Omega^c)}$.
This inequality follows from  Corollary 1 of Section 4.11.1 [M, p. 258].
Notice that in the notation of Maz'ja ${\cW}(\Omega^c) =
 W_{2, \frac{2(N-1)}{N-2}}(\Omega^c, \partial \Omega)$.   $\hfill //$

An immediate consequence of the above Lemma is that
the infimum in (\ref{2.2.2}) can be taken over $\cD^{1,2}(\Omega^c)$, that is
\begin{equation}
\mu_{\Omega}(\alpha) = 
\inf_{u\in   \cD^{1,2}(\Omega^c)}                
{\frac{\int_{\Omega^c}|\nabla u|^2dx -  \alpha\int_{\partial \Omega}   \frac{x \cdot \nu}{|x|^2}
  u^2 dS}{\int_{\Omega^c} \frac{u^2}{|x|^2} dx}}.
\end{equation}
We now state the analogue of Lemma \ref{lem:intf} for exterior domains.

\begin{lemma}
Let $\Omega$ be a bounded domain in $\R^N$, $N\geq 3$, containing the origin.  \nl
$\ia$ If $\alpha  \geq \frac{N-2}{2}$,  then $\mu_{\Omega}(\alpha)
 = \alpha(N-2-\alpha)$. Moreover, $|x|^{-\alpha} \in   \cD^{1,2}(\Omega^c)$
is a minimizer for  $\alpha >\frac{N-2}{2}$, whereas  for
$\alpha =\frac{N-2}{2}$ there is no $\cD^{1,2}(\Omega^c)$ minimizer. \nl
$\ib$ If  $0< \alpha < \frac{N-2}{2}$, and $\Omega$ starshaped with respect to zero, then 
$\mu_{\Omega}(\alpha ) = \Bigl(\frac{N-2}{2} \Bigr)^2$ and there is no $\cD^{1,2}(\Omega^c)$
  minimizer.
\label{lem:xe}
\end{lemma}       
{\em Proof:} The proof is quite  similar to the proof of the previous Lemmas \ref{lem:intf}
and \ref{lem:intf2}. An alternative proof can be given using the Kelvin transform; see 
the Remark that follows. $\hfill //$

\vspace{2mm}
\noindent {\bf Remark} There is a duality between the minimization problems
(\ref{2.2.1}) and (\ref{2.2.2}). Indeed, by means of the
Kelvin transform, $u(x) = |y|^{N-2} v(y)$, $y=x/|x|^2$, $x \in \Omega^c$  the domain
$\Omega^c$ is transformed
to a bounded domain containing the origin that we denote by $(\Omega^c)^*$. Denoting
by $\nu^*$ the outward pointing normal
to $\partial (\Omega^c)^*$ a straightforward calculation shows that
\[
\int_{\Omega^c} |\nabla_x u|^2 dx = \int_{(\Omega^c)^*}  |\nabla_y v|^2 dy
+(N-2) \int_{\partial (\Omega^c)^*} \frac{y \cdot \nu^*}{|y|^2} v^2 dS_y.
\]
Also,
\begin{eqnarray*}
\int_{\Omega^c} \frac{|u|^2}{|x|^2} dx & = & \int_{(\Omega^c)^*}  \frac{|v|^2}{|y|^2} dy, \\
\int_{\Omega^c} |u|^{\frac{2N}{N-2}} dx  & = & \int_{(\Omega^c)^*}  |v|^{\frac{2N}{N-2}} dy.
\end{eqnarray*} 
It can be seen  from these relations that $u \in {\cD^{1,2}(\Omega^c)}$ if and only
if $v \in  H^1((\Omega^c)^*)$. It then follows easily that
$
\mu_{\Omega}(\alpha) = \lambda_{(\Omega^c)^*}(N-2-\alpha)
$,
and that the existence  of  a minimizer for $\mu_{\Omega}(\alpha)$ in ${\cD^{1,2}(\Omega^c)}$ is
equivalent to the existence of a minimizer in $H^1(\Omega)$ for
$\lambda_{(\Omega^c)^*}(N-2-\alpha)$.

\subsection{Existence of minimizers}

In this section we establish a sufficient condition for the existence of
minimizers. We recall from Lemma \ref{lem:intf2} that when $\Omega$ is not starshaped
with respect to the origin,  $\alpha^*$ denotes the first zero of $\lambda_{\Omega}(\alpha)$.
We also set $\alpha^* = \infty$ in case $\Omega$ is starshaped with respect to zero.
Thus, in both cases we have $\lambda_{\Omega}(\alpha)>0$ for $0<\alpha<\alpha^*$.
Given $0< \alpha < \alpha^*$, and a nonnegative measurable potential $V$ we define
\begin{equation}
\lambda_{\Omega}(\alpha, V) :=
 \inf_{\scriptsize \begin{array}{c} 
   u \in H^1(\Omega) \\
   \int_{\Omega} V  u^2 dx >0
               \end{array}}
  \frac{\int_{\Omega}  |\nabla u|^2 dx + 
   \alpha\int_{\partial \Omega}   \frac{x \cdot \nu}{|x|^2}
  u^2 dS}{\int_{\Omega}  V u^2 dx}.
\label{concu}
\end{equation}
Note that with this notation $\lambda_{\Omega}(\alpha)=\lambda_{\Omega}(\alpha,|x|^{-2})$.
Since the numerator in (\ref{concu}) is always positive and finite when
$0<\alpha<\alpha^*$, we interpret $\lambda_{\Omega}(\alpha, V)=0$ in case there exists
$u\in H^1(\Omega)$ such that $\int_{\Omega} V  u^2 dx= +\infty$. It is worth mentioning
that $\lambda_{\Omega}(\alpha,V)$ is not monotone with respect to $\Omega$, unlike the
case of Dirichlet boundary conditions.

We denote by $B_r\subset\Omega$ the ball centered at zero with radius $r$.
We have the following
\begin{proposition}
Let $\Omega$ be a bounded domain in $\R^N$, $N\geq 3$, containing the origin, and let
$0 \leq V \in  L^{\frac{N}{2}}_{loc}(\overline{\Omega} \setminus\{0\})$. If 
for some $r>0$
\begin{equation}
0 <\lambda_{\Omega}(\alpha, V) < \lambda_{B_r}(\alpha, V)
\label{eq:2.14}
\end{equation}
then (\ref{concu}) has an $H^1(\Omega)$ minimizer.
\label{lem:2.3}
\end{proposition}
{\em Note.} It is a consequence of (\ref{eq:2.14}) that $\int_{\Omega}Vu^2<+\infty$ for
$u\in H^1(\Omega)$.

{\em Proof.} Let $\{u_j\} \in H^1(\Omega)$ be a minimizing sequence of the Rayleigh quotient
in (\ref{concu}).  We may normalize it so that $\int_{\Omega} V u_j^2 dx=1$. We claim that
$\|u_j\|_{ H^1(\Omega)} < C$. This will follow from two inequalities. The first
inequality follows
from the fact that $0<\alpha<\alpha^*$ and $\lambda_{\Omega}(\alpha^*)=0$ and reads
\begin{equation}
\int_{\Omega}  |\nabla u|^2 dx + 
   \alpha\int_{\partial \Omega}   \frac{x \cdot \nu}{|x|^2}
  u^2 dS \geq (1-\frac{\alpha}{\alpha^*})\int_{\Omega}  |\nabla u|^2 dx.
\label{2.10}
\end{equation}
The second one is a consequence of Lemma \ref{lem:intf} and reads
\begin{eqnarray}
\int_{\Omega}  |\nabla u|^2 dx + 
   \alpha\int_{\partial \Omega}   \frac{x \cdot \nu}{|x|^2}u^2 dS 
&\geq &  \lambda_{\Omega}(\alpha) \int_{\Omega} \frac{u^2}{|x|^2} dx \nonumber \\
 & \geq &K  \lambda_{\Omega}(\alpha) \int_{\Omega}u^2dx.
\label{2.21}
\end{eqnarray}
Thus, we may extract  a subsequence such that  $u_j \rightharpoonup u_0$
weakly in $ H^1(\Omega)$, and $u_j \rightarrow u_0$ strongly
in $L^p(\Omega)$, $1< p <\frac{2N}{N-2}$. Moreover, since $V\in L^{N/2}(\Omega\setminus B_r)$,
standard results give
\begin{equation}
\int_{\Omega\setminus B_r}Vu_j^2 dx \to \int_{\Omega\setminus B_r}Vu^2_0 dx.
\end{equation}
Also by the trace theorem
\begin{equation}
   \int_{\partial \Omega}   \frac{x \cdot \nu}{|x|^2}
  u^2_j dS  \rightarrow      \int_{\partial \Omega}   \frac{x \cdot \nu}{|x|^2}
  u^2_0 dS 
\label{2.22}
\end{equation}
Setting $u_j=v_j + u_0$ we easily see that as $j \rightarrow \infty$,
\begin{equation}
\int_{\Omega}  |\nabla u_j|^2 dx = \int_{\Omega}  |\nabla u_0|^2 dx +
 \int_{\Omega}  |\nabla v_j|^2 dx + o(1),
\label{2.23}
\end{equation}
and 
\begin{equation}
1= \int_{\Omega} V u_j^2 dx = \int_{\Omega} V u_0^2 dx +  \int_{\Omega} V v_j^2 dx + o(1).
\label{2.24}
\end{equation}
It then follows from (\ref{concu})  that
\begin{eqnarray}
\lambda_{\Omega}(\alpha, V)  &  =  &  \int_{\Omega}  |\nabla v_j|^2 dx +
\int_{\Omega}  |\nabla u_0|^2 dx + \alpha  \int_{\partial \Omega}   \frac{x \cdot \nu}{|x|^2}
  u^2_0 dS + o(1)  \nonumber   \\
& \geq &  \int_{\Omega}  |\nabla v_j|^2 dx + \lambda_{\Omega}(\alpha, V) \int_{\Omega} V u_0^2 dx
 + o(1).
\label{2.26}
\end{eqnarray}
We then have
\begin{eqnarray}
 \int_{\Omega}  |\nabla v_j|^2 dx &  \geq  &  \int_{B_r}  |\nabla v_j|^2 dx  \nonumber \\
& \geq &  \lambda_{B_r}(\alpha,V)
 \int_{B_r} V v_j^2 dx  \nonumber 
  -  \alpha  \int_{\partial B_r}   \frac{x \cdot \nu}{|x|^2}    v^2_j dS  \nonumber \\
& = &  \lambda_{B_r}(\alpha,V) \int_{B_r} V v_j^2 dx + o(1)   \nonumber \\
& = & \lambda_{B_r}(\alpha,V) \int_{\Omega} V v_j^2 dx + o(1) \nonumber \\
& = &  \lambda_{B_r}(\alpha,V) \Big(1- \int_{\Omega} V u_0^2 dx\Big) + o(1) 
\qquad (j\to\infty).
\end{eqnarray}
Using this and  (\ref{2.26}) we end up with
\begin{equation}
\Bigl( \lambda_{\Omega}(\alpha, V) - \lambda_{B_r}(\alpha, V)\Bigr)  
\Big(1- \int_{\Omega} V u_0^2 dx\Big) \geq 0,
\end{equation}
whence, since $\lambda_{\Omega}(\alpha,V)<\lambda_{B_r}(\alpha,V)$, it follows that
$\int_{\Omega} V u_0^2 dx \geq 1$. By lower semi continuity we conclude that
$\int_{\Omega} V u_0^2 dx = 1$. It then follows that $u_0$ is a minimizer for (\ref{concu}). 
$\hfill //$

As a consequence we have: \nl
{\em Completion of Proof of Lemma \ref{lem:intf2}(iii) (Existence of a minimizer)} :  Since 
$\frac{N-2}{2} \leq \alpha_* < \alpha < \alpha^*$
it follows from  Lemma \ref{lem:intf2} that $0<\lambda_{\Omega}(\alpha) < 
\bigl(\frac{N-2}{2} \big)^2 $.
If $B_r \subset \Omega$ is a ball centered at zero it follows by 
 Lemma \ref{lem:intf} that $\lambda_{B_r}(\alpha)=\bigl(\frac{N-2}{2} \big)^2 $. 
By Proposition \ref{lem:2.3},
$\lambda_{\Omega}(\alpha)$ is attained by an $H^1(\Omega)$ function.  $\hfill //$

We next state the corresponding result for the  exterior of a bounded  domain $\Omega$. 
For $0<a<N-2$ we define
\begin{equation}
\mu_{\Omega}(\alpha, V) :=
 \inf_{\scriptsize \begin{array}{c}
              u \in \cD^{1,2}(\Omega^c) \\
              \int_{\Omega^c} V  u^2 dx >0
               \end{array}}
  \frac{\int_{\Omega^c}  |\nabla u|^2 dx -
   \alpha\int_{\partial \Omega}   \frac{x \cdot \nu}{|x|^2}
  u^2 dS}{\int_{\Omega^c}  V u^2 dx}.
\label{cext}
\end{equation}

We then have
\begin{proposition}
Let $\Omega$ be a bounded domain in $\R^N$, $N\geq 3$, containing the origin, and let
$0 \leq V \in  L^{\frac{N}{2}}_{loc}(\overline{\Omega^c})$. 
Also let $0<\alpha< N-2$. If for some ball
$B_R \supset \Omega$  centered at zero 
\begin{equation}
0 <\mu_{\Omega}(\alpha, V) < \mu_{B_R}(\alpha, V),
\end{equation}
then (\ref{cext}) has a $\cD^{1,2}(\Omega^c)$ minimizer.
\label{lem:2.4}
\end{proposition}
The proof is similar to  that of the previous Proposition.



\setcounter{equation}{0}
\section{Hardy-Sobolev inequalities}
\subsection{Auxilliary inequalities}
We begin this section with two known Sobolev inequalities that will be used in the
sequel. In Theorems \ref{thm:gian} and \ref{thm:newsob} we then prove two new Sobolev
inequalities.

 By  the  classical inequality
of Caffarelli-Kohn-Nirenberg \cite{CKN} we have  that 
\begin{equation}
\int_{\R^N}|\nabla w|^2 |x|^{-2\alpha}dx \geq
c\Bigl(\int_{\R^N}|w|^{p}|x|^{-p\beta}dx
\Bigr)^{2/p},\quad w\in C^{\infty}_c(\R^N),
\label{eq:2.3}
\end{equation}
 with $p=2N/(N-2+2(\beta-\alpha))$,
provided $\alpha<(N-2)/2$ and $0\leq\beta -\alpha\leq 1$.

At the critical case $\alpha=\beta=(N-2)/2$ inequality (\ref{eq:2.3})
fails. A sharp substitute for bounded $\Omega$ was obtained in \cite{FT}, where it was
shown that, with
\[ X_1(t) = (1- \log t)^{-1} \;, ~~~~   t\in (0,1),\]
and $D=\sup_{\Omega}|x|$ there holds
\begin{equation}
\int_{\Omega}|\nabla w|^2|x|^{2-N}dx
\geq c\Bigl(\int_{\Omega}|w|^{\frac{2N}{N-2}}|x|^{-N}
X_1^{\frac{2N-2}{N-2}}(\frac{|x|}{D})dx\Bigr)^{(N-2)/N},
\label{eq:2.1}
\end{equation}
for all $w\in C^{\infty}_c(\Omega)$,
where the exponent $\frac{2N-2}{N-2}$ of $X_1(|x|/D)$ is optimal.

In the sequel  we will  make essential use of the following one-dimensional
result, which is a special case of a more general
statement by Maz'ja, cf. \cite[Theorem 3, p. 44]{M}:
\begin{proposition}
Let $A(r)$, $B(r)$ nonnegative functions such that $1/A(r)$  and $B(r)$
are integrable in $(r, \infty)$
and $(0,r)$ respectively, for all positive  $r < \infty$. Then, for  $q\geq 2$
the Sobolev inequality
\[  \int_0^{\infty} (v'(r))^2A(r)dr \geq c\Bigl( \int_0^{\infty} |v(r)|^q B(r)dr\Bigr)^{2/q}, \]
is valid for all $v\in C^{1}(0,\infty)$ that vanish near infinity,  if and only if
\[ \sup_{r>0}\; \Bigl(\int_0^r B(t)dt\Bigr) \; 
\Bigl(\int_r^{\infty}\frac{dt}{ A(t)} \Bigr)^{q/2} <+\infty. \]
\label{prop:maz}
\end{proposition}
The above proposition will be applied to higher dimensions by means of the
following
\begin{lemma}
Let $N \geq 2$. Suppose that $V \in L^{\infty}_{loc}(\R^N \setminus\{0\}) \cap L^1_{loc}(\R^N)$ 
is a radially symmetric function. We further  assume  that inequality
\begin{equation}
\int_{\R^N}|\nabla u|^2dx  -  \int_{\R^N}Vu^2dx \geq 0,
\label{eq:dia}
\end{equation}
is valid for all radially symmetric functions  $u\in C_c^{\infty}(\R^N)$. \nl
$\ia$ Then, (\ref{eq:dia}) is also  valid  for  nonradial functions, that is,
 for all  $u\in C_c^{\infty}(\R^N)$   \nl
$\ib$ If, in addition, 
\begin{equation}
0< \mbox{\rm ess \, sup}_{x \in \R^N} |x|^2 V(x) = \theta < \infty,
\label{2.11}
\end{equation}
then the following improved inequality holds
\begin{eqnarray} 
 & & \int_{\R^N}|\nabla u|^2dx  -\int_{\R^N}Vu^2dx  \geq \label{eq:dia2}  \\
& & ~~~ \geq  \int_{\R^N}|\nabla u_0|^2dx   - \int_{\R^N}Vu_0^2dx +
 \frac{N-1}{N-1+\theta} \,\,  \int_{\R^N}|\nabla(u-u_0)|^2 dx,   \nonumber
\end{eqnarray}
where $u_0(r)$ denotes the spherical average of $u\in C_c^{\infty}(\R^N)$, that is 
\begin{equation}
u_0(r) = \frac{1}{N \omega_N r^{N-1}} \int_{\partial B_r(0)} u(x)  dS_x.
\label{2.aver}
\end{equation}
\label{lem:2.2}
\end{lemma}
{\em Proof.}
Let $u\in C^{\infty}_c(\R^N)$ and let
\[u(x)=\sum_{m=0}^{\infty}f_m(\sigma)u_m(r),\]
be its decomposition into spherical harmonics; here $f_m$ are  orthogonal
in  $L^2(S^{N-1})$, normalized by
 $\frac{1}{N \omega_N} \int_{S^{N-1}} f_i(\sigma)f_j(\sigma) dS= \delta_{ij}$.
In particular $f_0(\sigma)=1$ and the first term in the above decomposition is given by 
(\ref{2.aver}).
The $f_m$'s are   eigenfunctions of the Laplace-Beltrami operator,
with corresponding eigenvalues
 $c_m=m(N-2+m)$, $m\geq 0$. An easy calculation shows that
$$
\int_{\R^N}(|\nabla u|^2-Vu^2)dx =
\sum_{m=0}^{\infty}\int_{\R^N}\Bigl\{|\nabla u_m|^2
+(\frac{c_m}{|x|^2}-V)u_m^2\Bigr\}dx=
$$
\begin{equation}
=\int_{\R^N}(|\nabla u_0|^2-Vu^2_0)dx + \sum_{m=1}^{\infty}\int_{\R^N}\Bigl\{|\nabla u_m|^2
+(\frac{c_m}{|x|^2}-V)u_m^2\Bigr\}dx.
\label{eq:au}
\end{equation}
Part (i) follows immediately since $c_m>0$ and $u_m$ are radially symmetric.

To prove part (ii) we first observe that
$$
 \int_{\R^N}|\nabla(u-u_0)|^2 dx =
\sum_{m=1}^{\infty}\int_{\R^N}
\Bigl\{|\nabla u_m|^2+\frac{c_m}{|x|^2}u_m^2\Bigr\}dx.
$$
In view of this and  (\ref{eq:au}) it is enough to establish that for any $m \geq 1$, there holds
\begin{equation}
\int_{\R^N}\Bigl\{|\nabla u_m|^2
+(\frac{c_m}{|x|^2}-V)u_m^2\Bigr\}dx \geq  \frac{N-1}{N-1+\theta}     \int_{\R^N}
\Bigl\{|\nabla u_m|^2+\frac{c_m}{|x|^2}u_m^2\Bigr\}dx .
\label{eq:zz}
\end{equation}
or, equivalently,
$$
\int_{\R^N}|\nabla u_m|^2dx \geq \int_{\R^N}  u_m^2   \Bigl\{ \frac{N-1+ \theta}{\theta} V -
 \frac{c_m}{|x|^2}  \Bigr\}  dx.
$$
Since $c_m \geq c_1=N-1$ it is enough to establish this for $c_m=N-1$.
By the definition of $\theta$,
cf (\ref{2.11}),
it follows easily that
$$
\frac{N-1+ \theta}{\theta} V -
 \frac{N-1}{|x|^2} \geq V,
$$
and the result follows from (\ref{eq:dia}).  $\hfill //$

As a consequence of this we next establish the following result.
\begin{lemma}
Let $N \geq 3$. Suppose that $V \in L^{\infty}_{loc}(\R^N \setminus\{0\}) \cap L^1_{loc}(\R^N)$ 
is a radially symmetric function, such that
$$
0< \mbox{\rm ess \,sup}_{x \in \R^N} |x|^2 V(x) = \theta < \infty,
$$
and $W \in L^{\infty}(\R^N)$ is a positive radially symmetric function. We further  assume that 
the inequality
\begin{equation}
\int_{\R^N}|\nabla u|^2dx  -  \int_{\R^N}Vu^2dx \geq
 c\left(\int_{\R^N}|u|^{2N/(N-2)}Wdx\right)^{(N-2)/N}
\label{eq:dia3}
\end{equation}
is valid for all radially symmetric functions  $u\in C_c^{\infty}(\R^N)$. Then 
inequality (\ref{eq:dia3}) is true for all $u\in C_c^{\infty}(\R^N)$ (without radial symmetry),
 provided the constant
 $c$  is  replaced by a new constant $C$
depending on $c$, $N$, $\theta$ and $\|W\|_{ L^{\infty}}$.
\label{lem:hol}
\end{lemma}

{\em Proof:} Starting from (\ref{eq:dia2}) we compute
\begin{eqnarray*}
&&\hspace{-2cm}\int_{\R^N}\left\{|\nabla u|^2 -Vu^2\right\}dx\geq \\
& \geq &\int_{\R^N}(|\nabla u_0|^2 -Vu_0^2)dx+  \frac{N-1}{N-1+\theta} \,\,
   \int_{\R^N}|\nabla(u-u_0)|^2dx\\
&\geq&c\left(\int_{\R^N}|u_0|^{2N/(N-2)}Wdx
\right)^{(N-2)/N}+ \\
&&+c' \left(\int_{\R^N}|u-u_0|^{2N/(N-2)}dx\right)^{(N-2)/N}\\
&\geq&C\left(\int_{\R^N}|u|^{2N/(N-2)}Wdx\right)^{(N-2)/N},
\end{eqnarray*}
where, for the last  inequalities we used the standard Sobolev inequality,  the boundedness
of $W$ and the triangle inequality.  $\hfill //$


\subsection{Hardy-Sobolev inequalities}
In this section we prove improved Hardy-Sobolev inequalities for
potentials that are critical either near zero
or near infinity. We first consider a potential which is critical near zero. 
 For $\epsilon>0$  we  define
\begin{equation}
V_{\epsilon}(x)=\darr{\Bigl(\frac{N-2}{2}\Bigr)^2|x|^{-2},}{|x|<1}
{\epsilon f(x),}{|x| \geq 1},
\label{eq:aa}
\end{equation}
where $f$ is  a non-negative, continuous and radially symmetric function on $\{|x| \geq 1\}$.
Moreover we assume $f$ to be subcritical, satisfying
\begin{equation}
f(x)\leq K|x|^{-2-\sigma} ,\qquad |x|\geq 1,  
\label{eq:subcr}
\end{equation}
for some $\sigma ,K>0$.

We  also define the following auxiliary function (cf. (\ref{eq:x}))
\begin{equation}
\tilde{X}_1(|x|)=\darr{X_1(|x|),}{|x|<1}{1,}{|x|>1}.
\label{eq:aa2}
\end{equation}

We shall henceforth denote by $B$ the unit ball in $\R^N$ centered at zero,
by $B^c$ its complement, and, as before,
we denote by $C^{\infty}_c(\R^N)|_{B^c}$ the set of restrictions on $B^c$ of all functions
$u\in C^{\infty}_c(\R^N)$.  We also  denote by $\nu$ the outward-pointing
unit vector on the surface $\partial B$. We have the following

\begin{theorem}
Let
\begin{equation}
\epsilon_0 =\inf_{u\in \cH^1(B^c)}
{\frac{\int_{B^c}|\nabla u|^2dx -\frac{N-2}{2}\int_{\partial B}u^2dS}{\int_{B^c}fu^2dx}}.
\label{eq:inf}
\end{equation}
Then $\epsilon_0>0$ and for any $\epsilon\in (0,\epsilon_0)$ there holds
\begin{equation}
\int_{\R^N}|\nabla u|^2dx -\int_{\R^N}V_{\epsilon}u^2dx 
\geq c\left(\int_{\R^N}|u|^{2N/(N-2)}\tilde{X}_1^{(2N-2)/(N-2)}(|x|)dx
\right)^{(N-2)/N},
\label{eq:nns}
\end{equation}
for all $u\in C^{\infty}_c(\R^N)$. Moreover, (\ref{eq:nns}) fails for $\epsilon=\epsilon_0$.
\label{thm:gian}
\end{theorem}
{\em Proof. } Since $f(x) \leq K |x|^{-2}$ the positivity of $\epsilon_0$ follows
from Lemma \ref{lem:xe} with $a=\frac{N-2}{2}$, yielding in fact
$\epsilon_0 \geq K^{-1}((N-2)/2)^2.$

Let us now fix $\epsilon \in (0,\epsilon_0)$. 
By Lemma \ref{lem:hol}, it is enough to prove (\ref{eq:nns}) in the case where
$u$ is radially symmetric, $u=u(r)$.
Now, there exists a radially symmetric
and positive function $\tilde{\psi}$ on $B^c$ which solves the Robin problem
\[\darrn{\Delta\tilde{\psi} +\epsilon f\tilde{\psi}=0,}{|x|>1,}
{\parder{\tilde{\psi}}{\nu}=-\frac{N-2}{2}\tilde{\psi},}{|x|=1.}\]
The existence of such a $\tilde{\psi}$ can be easily derived, for example,
by a shooting argument from $\{|x|=1\}$.
We assume that $\tilde{\psi}$ is normalized so that $\tilde{\psi}=1$ on $\{|x|=1\}$.
The function
\begin{equation}
\psi(x)=\darr{|x|^{-(N-2)/2},}{|x|<1,}{\tilde{\psi}(x),}{|x|>1,}
\label{eq:psi}
\end{equation}
then lies in $C^1(\R^N\setminus\{0\})$, is positive,
radially symmetric and satisfies $\Delta\psi+V_{\epsilon}\psi=0$ in $\R^N$.
Following \cite{FT} we change variables, $u=\psi v$,
and (\ref{eq:nns}) for radially symmetric functions is then written as
\begin{equation}
\int_0^{\infty}(v')^2\psi^2r^{N-1}dr\geq c\left(\int_0^{\infty}|v|^{2N/(N-2)}
\psi^{2N/(N-2)}\tilde{X}_1^{2(N-1)/(N-2)}dr\right)^{(N-2)/N}.
\label{eq:papa}
\end{equation}
We claim that $\psi(r)$ has a positive limit as $r\to +\infty$. Indeed,
since $(r^{N-1}\psi')'=-r^{N-1}V_{\epsilon}\psi<0$ and $\psi'(1)<0$, $\psi(r)$
is decreasing on $(1,+\infty)$. 
If the limit $\lim_{r\to +\infty}\psi(r)$ were zero it would then
follow from \cite[Theorem 2.9]{LN} that $\psi(r) < cr^{2-N}$ near infinity,
which then easily implies $\psi\in\cH^1(B^c)$.
Hence $\psi$ can be taken as a test function for the infimum
in the right-hand side of (\ref{eq:inf}), in which case the value
of the Rayleigh quotient is $\epsilon<\epsilon_0$, contradicting
the definition of $\epsilon_0$. Hence $\lim_{r\to +\infty}\psi(r)=l>0$.
Using this we deduce (\ref{eq:papa}) from Proposition \ref{prop:maz}.
Hence (\ref{eq:nns}) has been proved.

We finally show that (\ref{eq:nns}) fails for $\epsilon=\epsilon_0$. For this we will
use Proposition \ref{lem:2.4}. Let $B_R\supset B_1$. Then
$V_{\epsilon_0}(x) \leq \epsilon_0   KR^{-\sigma}|x|^{-2}$ on $(B_R)^c$, hence
\[ \mu_{B_R}(\alpha,V_{\epsilon_0})\geq\frac{R^{\sigma}}{K}\mu_{B_R}(\alpha).\]
By Lemma \ref{lem:xe} (i), for  $\alpha \in (0,\, N-2)$
 $\mu_{B_R}(\alpha)$, 
is positive and independent
of $R$;  taking $R$ large enough we
 have $\mu_B(\alpha,V_{\epsilon_0})<\mu_{B_R}(\alpha,V_{\epsilon_0})$.
Hence, by Proposition \ref{lem:2.4} -- with $\alpha=(N-2)/2$ --
there exists an  $\cH^1(B^c)$-minimizer $\phi$ to (\ref{eq:inf}).
It is standard to show that
$\phi$ is simple, radial and of one sign; we normalize it by $\phi|_{|x|=1}=1$
and for $\theta>0$ we define the function $u_{\theta}\in H^1(\R^N)$ by
\[u_{\theta}(x)=\darr{|x|^{-\frac{N-2}{2}+\theta},}{|x|<1;}{\phi(x),}{|x|>1.}\]
We then compute the left-hand side of (\ref{eq:nns}): in $B$ there holds
$\Delta u_{\theta}+V_{\epsilon_0}u_{\theta}
=\theta^2 u_{\theta}$, hence
\begin{eqnarray*}
&&\int_{\R^N}(|\nabla u_{\theta}|^2 -V_{\epsilon_0}u_{\theta}^2)dx =\\
&=&-\int_B(u_{\theta}\Delta u_{\theta}+V_{\epsilon_0}u_{\theta}^2)dx+
\int_{\partial  B}u_{\theta}\parder{u_{\theta}}{\nu}dS +\int_{B^c}(|\nabla u_{\theta}|^2 
-V_{\epsilon_0}u_{\theta}^2)dx \\
&=&-\theta^2\int_B\frac{u_{\theta}^2}{r^2}dx+N\omega_N (-\frac{N-2}{2}+\theta)
+\frac{N-2}{2}\int_{\partial B}u_{\theta}^2dS\\
&=&\frac{N\omega_N\theta}{2}.
\end{eqnarray*}
On the other hand for the right-hand side of  (\ref{eq:nns}) we have
\[\int_{\R^N}u_{\theta}^{2N/(N-2)}\tilde{X}_1^{\frac{2N-2}{N-2}}dx\geq
\int_{B^c}u_{\theta}^{2N/(N-2)}dx
=\int_{B^c}\phi^{2N/(N-2)},\]
the last term being independent of $\theta$.
Letting $\theta\to 0$ we conclude that (\ref{eq:nns}) fails for $\epsilon=\epsilon_0$.
$\hfill //$

We close this section proving a Sobolev inequality which involves radial potentials
with critical behaviour $((N-2)/2)^2|x|^{-2}$ near infinity.
Let $g$ be  a non-negative, radially symmetric, and  continuous in $B \setminus \{0\}$
function that  is subcritical near zero, satisfying
\[ g(x)\leq K|x|^{-2+\sigma} ,\qquad |x|\leq 1,\]
for some $\sigma , \, K>0$. For $\epsilon>0$ we define
\begin{equation}
\hat{V}_{\epsilon}(x)=\darr{\epsilon g(x),}{|x|<1,}{(\frac{N-2}{2})^2|x|^{-2},}{|x|>1.}
\label{eq:xg}
\end{equation}
We also set
\[\tilde{Y}_1(|x|)= \darr{1,}{|x|<1,}{(1+\ln |x|)^{-1},}{|x|>1.}\]
We then have
\begin{theorem}
Let
\[\bar{\epsilon}_0=\inf_{u\in H^1(B)}
{\frac{\int_{B}|\nabla u|^2dx +\frac{N-2}{2}\int_{\partial B}u^2dS}{\int_{B}gu^2dx}}. \]
Then $\bar{\epsilon}_0>0$ and for any $\epsilon\in (0,\bar{\epsilon}_0)$ there holds
\begin{equation}
\int_{\R^N}|\nabla u|^2dx -\int_{\R^N}\hat{V}_{\epsilon}u^2dx 
\geq c\left(\int_{\R^N}|u|^{2N/(N-2)}\tilde{Y}_1^{2(N-1)/(N-2)}(|x|)dx
\right)^{(N-2)/N},
\label{eq:newsob}
\end{equation}
for all $u\in C^{\infty}_c(\R^N)$.
Moreover, (\ref{eq:newsob}) fails for $\epsilon=\bar{\epsilon}_0$.
\label{thm:newsob}
\end{theorem}
{\em Proof. }The proof of (\ref{eq:newsob}) follows closely that of
Theorem \ref{thm:gian}, reversing essentially
the role of $B$ and $B^c$ while making the necessary adjustments;
in particular, we now use Lemma \ref{lem:intf} instead of
Lemma \ref{lem:xe}. In fact, an alternative
and simpler proof consists in simply taking the Kelvin transform of (\ref{eq:nns}). The
optimality of $\bar{\epsilon}_0$ is also proven analogously; we omit the details. $\hfill //$




\setcounter{equation}{0}
\section{Heat kernel estimates}

In this section we shall apply the Sobolev inequalities of Section 3 to obtain
heat kernel estimates for the Schr\"{o}dinger operator
\[Hu=-\Delta u- Vu,\quad u|_{\partial\Omega}=0,\]
for various critical potentials $V$. We always assume that $\Omega$ is a domain
containing the origin in $\R^N$, $N\geq 3$. 
We shall consider the case of bounded $\Omega$, as well as
 the case $\Omega=\R^N$.
The operator $H$ is defined via quadratic forms, with initial domain
$C^1_c(\Omega\setminus\{0\})$; it will always be the case that $H\geq 0$.
Note that, equivalently, we could have set $C^1_c(\Omega)$ as the initial domain.

We shall use the standard technique of transference to a
weighted $L^2$ space, which we now describe briefly.
Let $\phi\in C^1(\Omega\setminus\{0\})$ be positive and 
such that $\Delta\phi\in L^1_{loc}(\Omega\setminus\{0\})$.
The unitary map
\begin{equation}
L^2 (\Omega)\ni u\mapsto w= \frac{u}{\phi} \in L^2_{\phi}:=L^2(\Omega,\phi^2\, dx)
\label{eq:br}
\end{equation}
satisfies
\begin{equation}
\int_{\Omega}(|\nabla u|^2 - Vu^2)dx
=\int_{\Omega}\Bigl(|\nabla w|^2 -Vw^2-\frac{\Delta\phi}{\phi}w^2\Bigr)\phi^2dx
\label{eq:20a}
\end{equation}
for all $u\in C^1_c(\Omega\setminus\{0\})$. Hence, if in addition $\phi$ satisfies
$\Delta\phi+V\phi=0$ (weakly) on $\Omega\setminus\{0\}$, then
\begin{equation}
\int_{\Omega}(|\nabla u|^2 - Vu^2)dx
=\int_{\Omega}|\nabla w|^2\phi^2dx
\label{eq:20b}
\end{equation}
for all $u\in C^1_c(\Omega\setminus\{0\})$. Hence $H$ is unitarily equivalent
via (\ref{eq:br}) to the self-adjoint operator $H_{\phi}$ on $L^2_{\phi}$,
defined initially on $C^1_c(\Omega\setminus\{0\})$ and given formally by
\[ H_{\phi}w=-\frac{1}{\phi^2}\diver(\phi^2\nabla w),\quad w|_{\partial\Omega}=0 .\]
The space $C^1_c(\Omega\setminus\{0\})$ is invariant under
multiplication (or division) by $\phi$ and hence it is a form core also for $H_{\phi}$.
Moreover, a Sobolev inequality of the form
\[\int_{\Omega}(|\nabla u|^2 - Vu^2)dx \geq c \left(\int_{\Omega}|u|^qWdx\right)^{2/q}\]
is valid for all $u\in C_c^1(\Omega\setminus\{0\})$ if and only if
\[\int_{\Omega}|\nabla w|^2\phi^2dx \geq c  \left(\int_{\Omega}|w|^q\phi^qWdx\right)^{2/q}\]
for all $w\in C_c^1(\Omega\setminus\{0\})$.
Finally, the heat kernels of $H$ and $H_{\phi}$ are related by
\begin{equation}
K(t,x,y)=\phi(x)\phi(y)K_{\phi}(t,x,y),\quad 
t>0,\; x,y\in\Omega ,
\label{eq:2.20}
\end{equation}
and hence one can obtain estimates on $K(t,x,y)$ via estimates on $K_{\phi}(t,x,y)$.

\vspace{10pt}
\noindent
{\bf Example.} As a typical example let us consider the case of a bounded domain $\Omega$
in $\R^N$, $N\geq 3$, and let $V(x)=\lambda|x|^{-2}$, $\lambda\leq ((N-2)/2)^2$.
Let $\phi(x)=|x|^{-\alpha}$, $\alpha$ being the smallest solution of
$\alpha(N-2-\alpha)=\lambda$. Then $\Delta \phi+V\phi=0$ on $\Omega\setminus\{0\}$
and therefore
\begin{equation}
\int_{\Omega}(|\nabla u|^2 -\lambda\frac{u^2}{|x|^2})dx
=\int_{\Omega}|\nabla w|^2|x|^{-2\alpha}dx.
\label{eq:2a}
\end{equation}
for all $u\in C^1_c(\Omega\setminus\{0\})$ or, equivalently,
for all $w\in C^1_c(\Omega\setminus\{0\})$.
Moreover, the heat kernel of $H=-\Delta-V$ is related to the heat kernel of $H_{\phi}$ by
\[ K(t,x,y)=|x|^{-\alpha}|y|^{-\alpha}K_{\phi}(t,x,y).\]
We note here that a simple approximation argument shows that for $\lambda<((N-2)/2)^2$
the form domain of $H$ is $H^1_0(\Omega)$,
but at the critical case $\lambda=((N-2)/2)^2$ the form domain
is strictly larger than $H^1_0(\Omega)$; see also \cite{FT}.

Sobolev inequalities are related to heat kernel estimates by the following standard
result \cite[Theorem 2.4.2]{D}: for any $q>2$,

\begin{equation}
\left\{
\begin{array}{l}
\mbox{the upper bound} \\
\qquad\qquad K_{\phi}(t,x,y)< ct^{-q/2},\quad t>0,\, x,y\in\Omega \\
\mbox{is equivalent to the Sobolev inequality} \\
\qquad\qquad \int_{\Omega}|\nabla w|^2\phi^2dx\geq c\Bigl(\int_{\Omega}
|w|^{\frac{2q}{q-2}}\phi^2 dx\Bigr)^{(q-2)/q},\;\; w\in C^1_c(\Omega\setminus\{0\}).
\end{array}\right.
\label{eq:eq}
\end{equation}

In the rest of this
section we shall apply the Hardy-Sobolev inequalities of Section 3 in order
to obtain upper estimates on the heat kernel $K(t,x,y)$ of the operator $-\Delta -V$ for
critical and subcritical potentials $V$. For this we shall use (\ref{eq:2.20}) for appropriate
functions $\phi$, together with uniform estimates on $K_{\phi}(t,x,y)$, obtained by means of
(\ref{eq:eq}). We present on-diagonal estimates, in the sense that we do not include the usual
Gaussian term, which can be added by standard methods.

We assume that $\Omega$ is a domain in $\R^N$, $N\geq 3$. We retain the notation
introduced in the last example, and, in particular, we have $H=-\Delta-\lambda |x|^{-2}$,
subject to Dirichlet boundary conditions on $\partial\Omega$. Part $\ia$ of the next
proposition has been proved in \cite{MS} but we include a proof for the sake of
completeness. Part $\ib$ improves upon the corresponding estimate of \cite{MS}, who had
$K(t,x,y)<ct^{-\frac{N-\alpha}{2}-\epsilon}|x|^{-\alpha}|y|^{-\alpha}$.

\begin{proposition}{\bf(subcritical case)} Let  $K(t,x,y)$ be  the heat kernel of
  $H=-\Delta-\lambda \frac{1}{|x|^{2}}$,
subject to Dirichlet boundary conditions on $\partial\Omega$. For 
$\lambda<((N-2)/2)^2$, let $\alpha$ be the smallest solution of
 $\alpha(N-2-\alpha)=\lambda$.   \nl
$\ia$ If $\Omega$ is bounded and $0\leq\lambda<((N-2)/2)^2$ then
\begin{equation}
K(t,x,y)< ct^{-N/2}|x|^{-\alpha}|y|^{-\alpha},\quad
t>0,\, x,y\in\Omega.
\label{eq:est1}
\end{equation}
$\ib$ For any $\Omega \subset \R^N$ (bounded or unbounded) and  $\lambda\leq 0$ then
\begin{equation}
K(t,x,y)< ct^{-N/2}\min\{1, (\frac{|x|}{t^{1/2}})^{-\alpha}\}
\min\{1, (\frac{|y|}{t^{1/2}})^{-\alpha}\},\quad
t>0,\, x,y\in\Omega.
\label{eq:est2}
\end{equation}
\label{prop:pos}
\end{proposition}
\vspace{-.5cm}
{\em Proof.} (i) The boundedness of $\Omega$ together with
(\ref{eq:2.3}) imply
\begin{equation}
\int_{\Omega}|\nabla w|^2 |x|^{-2\alpha}dx \geq
c\Bigl(\int_{\Omega}|w|^{\frac{2N}{N-2}}|x|^{-2\alpha}dx
\Bigr)^{(N-2)/N},\quad w\in C^{\infty}_c(\Omega).
\label{eq:2.2a}
\end{equation}
By (\ref{eq:eq}) this implies $K_{\phi}(t,x,y)<ct^{-N/2}$, from which
(\ref{eq:est1}) follows using (\ref{eq:2.20}).

(ii) Comparison with the Laplacian implies that $K(t,x,y)<ct^{-N/2}$.
Moreover, inequality (\ref{eq:2.3}) for $\beta p=2\alpha$ reads
\[\int_{\Omega}|\nabla w|^2|x|^{-2\alpha}dx\geq
c\Bigl(\int_{\Omega}|w|^{\frac{2(N-2\alpha)}{N-2-2\alpha}}|x|^{-2\alpha} 
dx\Bigr)^{\frac{N-2-2\alpha}{N-2\alpha}}.\]
By means of (\ref{eq:eq}) we deduce that  
 $K_{\phi}(t,x,y)< ct^{-N/2+\alpha}$, $t>0,\, x,y\in\Omega$. Hence
\begin{equation}
K(t,x,y)< ct^{-N/2}\min\{1, (\frac{|x|}{t^{1/2}})^{-\alpha}
(\frac{|y|}{t^{1/2}})^{-\alpha}\},\quad
t>0,\; x,y\in\Omega .
\label{eq:est3}
\end{equation}
This proves (\ref{eq:est2}) when $x=y$. The general case follows from
the semigroup property since
\begin{eqnarray*}
K(t,x,y)&=&\int_{\Omega}K(t/2,x,z)K(t/2,z,y)dz \\
&\leq& \Bigl(\int_{\Omega}K(t/2,x,z)^2dz\Bigr)^{1/2}
\Bigl(\int_{\Omega}K(t/2,z,y)^2dz\Bigr)^{1/2} \\
&=&K(t,x,x)^{1/2}K(t,y,y)^{1/2}.
\end{eqnarray*}
$\hfill //$


\begin{theorem}
{\bf (critical case)} Let  $\Omega$ be a  bounded domain and 
 $K(t,x,y)$ be  the heat kernel of
  $H=-\Delta- (\frac{N-2}{2})^2  \frac{1}{|x|^{2}}$,
subject to Dirichlet boundary conditions on $\partial\Omega$. Then 
\begin{equation}
K(t,x,y) <ct^{-\frac{N}{2}}|x|^{-\frac{N-2}{2}}|y|^{-\frac{N-2}{2}}, \quad t>0,\, x,y\in\Omega.
\label{eq:4.4}
\end{equation}
\label{thm:mainthm}
\end{theorem}
\vspace{-.5cm}
{\em Proof. }Estimate (\ref{eq:2.1}) implies the weaker inequality
\begin{equation}
\int_{\Omega}|\nabla w|^2|x|^{2-N}dx
\geq c\Bigl(\int_{\Omega}|w|^{\frac{2N}{N-2}}|x|^{2-N}
dx\Bigr)^{(N-2)/N}.
\label{eq:4.2}
\end{equation}
Hence $K_{\phi}(t,x,y) <ct^{-\frac{N}{2}}$ and (\ref{eq:4.4}) follows.$\hfill //$

We next consider the case $\Omega = \R^N$, and the potential is critical at zero. 
More precisely we consider the potential  $V_{\epsilon}$ defined by (\ref{eq:aa}),
that is:
\[
V_{\epsilon}(x)=\darr{\Bigl(\frac{N-2}{2}\Bigr)^2|x|^{-2},}{|x|<1}
{\epsilon f(x),}{|x| \geq 1},
\]
where $f$ is  a non-negative, continuous and radially symmetric function on $\{|x| \geq 1\}$.
Moreover we assume $f$ to be subcritical, that is it satisfies (\ref{eq:subcr}):
\[f(x)\leq K|x|^{-2-\sigma} ,\qquad |x|\geq 1,\]
for some $\sigma, K>0$.


 We retain the notation  of section  3.1, and in particular we recall the definition
(\ref{eq:inf}) of $\epsilon_0$. We have
\begin{theorem}
{\bf (the operator $-\Delta-V_{\epsilon}$ on $\R^N$)}
For any $\epsilon\in (0,\epsilon_0)$ the heat kernel of the operator $-\Delta -V_{\epsilon}$
satisfies
\begin{equation}
K(t,x,y)<ct^{-N/2}\max\{|x|^{-\frac{N-2}{2}},1\}\max\{|y|^{-\frac{N-2}{2}},1\}
,\quad t>0,\; x,y\in\R^N.
\label{eq:elec}
\end{equation}
\label{thm:elec}
\end{theorem}
{\em Proof.} Let $\psi(x)$ be as in the proof of Theorem \ref{thm:gian}, cf (\ref{eq:psi}).
It follows from (\ref{eq:nns}) that
\begin{equation}
\int_{\R^N}|\nabla v|^2\psi^2dx\geq c\left(\int_{\R^N}|v|^{2N/(N-2)}\psi^{2N/(N-2)}
\tilde{X}_1^{(2N-2)/(N-2)}dx\right)^{(N-2)/N},
\label{eq:hal}
\end{equation}
for all $v\in C^1_c(\R^N\setminus\{0\})$. Since
$\psi^{2N/(N-2)}\tilde{X}_1^{(2N-2)/(N-2)} \geq c\psi^2$
and $C^1_c(\R^N\setminus\{0\})$ is a form core for $H_{\psi}$ and 
we conclude that $K_{\psi}(t,x,y)<ct^{-N/2}$ whence,
\[
K(t,x,y) <ct^{-N/2}  {\psi}(x) {\psi}(y).
\]
The required estimate on $K(t,x,y)$ follows if we note that
\[ c_1 \max\{|x|^{-(N-2)/2},1\}\leq \psi(x) \leq c_2 \max\{|x|^{-(N-2)/2},1\}\;\; ,
\quad x\in\R^N.\]
$\hfill //$

It is well known that the estimates of the above theorems can be
improved to yield Gaussian decay of the heat kernel. We have
\begin{proposition}
Proposition \ref{prop:pos} as well as Theorems \ref{thm:mainthm} and \ref{thm:elec}
can be improved by adding a factor $c_{\delta}\exp\{-|x-y|^2/((4+\delta)t)\}$
to the corresponding right-hand sides.
\label{prop:bern}
\end{proposition}
{\em Proof. }The proof is standard. One can use Davies's method
of exponential perturbation \cite{D} or Theorem 1.1 of
\cite{G}. Note that the argument is applied to the operator $H_{\phi}$ -- not
to $H$. $\hfill //$

\section{Logarithmic refinements}

Our aim in this section is to obtain refined versions of the
improved Hardy-Sobolev inequalities of Section 3. As an application,
we prove heat kernel estimates for $H=-\Delta -V-V_1$ where
$V$ is one of the potentials studied in Section 4 (that is,  $V(x)=((N-2)/2)^2|x|^{-2}$ near
zero) but $V_1>0$ is also critical. 
The criticality of $V_1$ is meant in the sense that the following 
improved Hardy inequality holds
\[ \int_{\Omega}|\nabla u|^2 -  \int_{\Omega}Vu^2 
\geq
\int_{\Omega}V_1u^2,
\qquad u\in C^{\infty}_c(\Omega),\]
whereas this inequality is no longer true if we replace $V_1$ by $(1+\epsilon)V_1$.
Of course, $V_1$ is of lower order with respect to $|x|^{-2}$ (near $x=0$) since $((N-2)/2)^2|x|^{-2}$ is
 already critical
for the validity of the (simple) Hardy inequality. It is remarkable that the addition of the extra potential $V_1$
does not affect the time dependence of the heat kernel estimates, but only affects the spatial singularity
at the origin, which is increased by a logarithmic factor; see Theorems \ref{thm:lognon} and \ref{thm:hk100}.

More precisely, recalling that $X_1(t) = (1- \log t)^{-1}$, let us introduce the functions
\begin{equation}
X_{k+1}(t) = X_{1}(X_k(t)), \quad k=1,2,\ldots,\quad  t\in (0,1).
\label{eq:x}
\end{equation}
These are iterated logarithmic functions that vanish at an increasingly slow rate at $t=0$
and are equal to one at $t=1$. In \cite{FT} the following improved Hardy inequality was
obtained for a bounded domain $\Omega$ with $D=\sup_{\Omega}|x|$ :
\begin{eqnarray}
&& \int_{\Omega} \Bigl\{ |\nabla u|^2 -\Bigl(\frac{N-2}{2}\Bigr)^2\frac{u^2}{|x|^2}
-\frac{u^2}{4|x|^2}\sum_{i=1}^k 
X_1^2(\frac{|x|}{D})\ldots X_i^2(\frac{|x|}{D})\Bigr\}dx \nonumber \\
&& \hspace{2cm}\geq \frac{1}{4}\int_{\Omega}\frac{u^2}{|x|^2}X_1^2(\frac{|x|}{D})
\ldots X_{k+1}^2(\frac{|x|}{D})dx,
\qquad\quad u\in C^{\infty}_c(\Omega).
\label{eq:5.01}
\end{eqnarray}
The potentials in the left-hand side of (\ref{eq:5.01}) are critical for each $k$, in the
sense that (\ref{eq:5.01}) is sharp:
the term $X_{k+1}^2$ cannot be replaced by $c_{\epsilon}X_{k+1}^{2-\epsilon}$ for any $\epsilon>0$, and the constant
$1/4$ in the right-hand side is also optimal. In Theorem \ref{thm:lognon}
and for bounded $\Omega$ we obtain upper estimates on the heat kernel of the operator
\begin{equation}
H=-\Delta -\Bigl(\frac{N-2}{2}\Bigr)^2
\frac{1}{|x|^2}-\frac{1}{4|x|^2}\sum_{i=1}^{k-1}
X_1^2\ldots X_i^2  -\frac{\mu}{|x|^2}X_1^2\ldots X_k^2
\label{eq:anti}
\end{equation}
for $\mu> 1/4$, as well as for the critical case $\mu=1/4$; for this we use
results obtained in \cite{FT}.
For the critical case $\mu=1/4$ we also consider operators defined on $\R^N$,
in analogy to the operator $-\Delta-V_{\epsilon}$ of Theorem \ref{thm:elec};
for this we use Theorem \ref{thm:gian1} below, and the corresponding heat kernel
estimate is given in Theorem \ref{thm:hk100}.

\subsection{Refined Hardy-Sobolev inequalities}

In this subsection we prove two theorems that are refined
versions of Theorems \ref{thm:gian} and \ref{thm:newsob} correspondingly. We recall
definition (\ref{eq:x}) and set 
\[\tilde{X}_k(|x|)=\darr{X_k(|x|),}{|x|<1,}{1,}{|x|>1,}\]
\[ \darr{Y_k(|x|)=X_k(1/|x|),}{|x|>1,}{\tilde{Y}_k(|x|)=\tilde{X}_k(1/|x|),}{|x|> 0.}\]
We point out the differentiation rules for $X_k(r)$ and $Y_k(r)$:
\begin{equation}
\frac{d}{dr}X_k^a=\frac{a}{r}X_1\ldots X_{k-1}X_k^{a+1},\quad
\frac{d}{dr}Y_k^a=-\frac{a}{r}Y_1Y_2\ldots Y_{k-1}Y_k^{a+1}, \quad \quad r=|x|,
\label{eq:dr}
\end{equation}
valid for $0<r<1$ and $r>1$ respectively, which are easily proved by induction.

As in Theorem \ref{thm:gian}, we assume that
$f$ is a non-negative, continuous and radially symmetric function on $B^c$ satisfying
(\ref{eq:subcr}), that is,
\[ f(x)\leq K|x|^{-2-\sigma} ,\qquad |x|\geq 1,\]
for some $\sigma ,K>0$. For $\epsilon>0$ we also define
\begin{equation}
V_{k,\epsilon}(x)=\darr{\Bigl(\frac{N-2}{2}\Bigr)^2|x|^{-2}+\frac{1}{4}|x|^{-2}
\sum_{i=1}^kX_1^2(|x|)\ldots X_i^2(|x|),}{|x|<1,}
{\epsilon f(x),}{|x|>1.}
\label{eq:vke}
\end{equation}
We then have
\begin{theorem}
Assume that $k<N-2$ and define
\begin{equation}
\epsilon_{k,0} =\inf_{u\in \cH^1(B^c)}
{\frac{\int_{B^c}|\nabla u|^2dx -\frac{N-2+k}{2}\int_{\partial B}u^2dS}
{\int_{B^c}fu^2dx}}.
\label{eq:inf1}
\end{equation}
Then $\epsilon_{k,0}>0$ and for $\epsilon\in (0,\epsilon_{k,0})$ there holds
\begin{equation}
\int_{\R^N}|\nabla u|^2dx -\int_{\R^N}V_{k,\epsilon}u^2dx 
\geq c\left(\int_{\R^N}|u|^{2N/(N-2)}(\tilde{X}_1\ldots
\tilde{X}_{k+1})^{\frac{2N-2}{N-2}}dx
\right)^{(N-2)/N}\hspace{-.2cm},
\label{eq:nns1}
\end{equation}
for all $u\in C^{\infty}_c(\R^N)$.
\label{thm:gian1}
\end{theorem}
{\bf  Remark.} The constant $\epsilon_{k,0}$ is optimal in the sense that inequality
(\ref{eq:nns1}) fails for $\epsilon = \epsilon_{k,0}$. Also  the exponent $(2N-2)/(N-2)$
in (\ref{eq:nns1})
is sharp in the sense that it cannot be replaced by a smaller exponent.
 The proof of
these two facts is  rather involved; see [FT] for similar arguments. We do not use these
 facts in the sequel. \nl
{\em Proof.} The proof follows closely that of Theorem \ref{thm:gian}, so we
only give a sketch of it. The positivity of $\epsilon_{k,0}$ follows from
Lemma \ref{lem:xe} (i), yielding $\epsilon_{k,0}\geq K^{-1}\mu_B((N-2+k)/2)$.
Now let $\epsilon\in (0,\epsilon_{k,0})$ be fixed and let $\tilde{\psi}$ be the
radially symmetric solution to the problem
\[\darrn{\Delta\tilde{\psi} +V_{k,\epsilon} \tilde{\psi}=0,}{|x|>1,}
{\parder{\tilde{\psi}}{\nu}=-\frac{N-2+k}{2}\tilde{\psi},}{|x|=1,}\]
normalized so that $\tilde{\psi}=1$ on $\{|x|=1\}$.
The function
\begin{equation}
\psi(x)=\darr{|x|^{-(N-2)/2}X_1^{-1/2}\ldots X_k^{-1/2},}{|x|<1}{\tilde{\psi}(x),}{|x|>1,}
\label{eq:defpsi}
\end{equation}
is then $C^1$, radially symmetric and a direct computation
which uses (\ref{eq:dr}) shows that $\Delta\psi+V_{k,\epsilon}\psi=0$ in $\R^N$.
Exactly as in Theorem \ref{thm:gian}, $\psi$ is positive, radially symmetric and
has a positive limit as $r\to+\infty$.
We then prove (\ref{eq:nns1}) in the case where $u$ is radially symmetric,
using once again Proposition \ref{prop:maz}. The validity of (\ref{eq:nns1})
for general $u\in C^{\infty}_c(\R^N)$ follows from Lemma \ref{lem:hol}. $\hfill //$

We finally prove a refined version of Theorem \ref{thm:newsob}.
Let us fix a a non-negative, continuous and radially symmetric
function  $g$ on $B=\{|x|<1\}$, such that
\[g(x)\leq K|x|^{-2+\sigma} ,\qquad |x|< 1,\]
for some $\sigma ,K>0$. Further for $\epsilon>0$ we define
\[\hat{V}_{k,\epsilon}(x)=\darr{\epsilon g(x),}{|x|<1,}{\Bigl(\frac{N-2}{2}\Bigr)^2|x|^{-2}
+\frac{1}{4|x|^2}\sum_{i=1}^kY_1^2(|x|)\ldots Y_i^2(|x|),}{|x|>1.}\]

We then have
\begin{theorem}
Assume that $k<N-2$ and define
\[\bar{\epsilon}_{k,0}=\inf_{u\in H^1(B)}
{\frac{\int_{B}|\nabla u|^2dx +\frac{N-2-k}{2}\int_{\partial B}u^2dS}{\int_{B}gu^2dx}}. \]
Then $\bar{\epsilon}_{k,0}>0$ and for $\epsilon\in (0,\bar{\epsilon}_{k,0})$ there holds
\begin{equation}
 \int_{\R^N}|\nabla u|^2dx -\int_{\R^N}\hat{V}_{k,\epsilon}u^2dx 
\geq  c\left(\int_{\R^N}|u|^{2N/(N-2)}(\tilde{Y}_1\ldots
\tilde{Y}_{k+1})^ {\frac{2N-2}{N-2}}  dx
\right)^{(N-2)/N},
\label{eq:newsob1}
\end{equation}
for all $u\in C^{\infty}_c(\R^N)$.
\label{thm:newsob1}
\end{theorem}
{\em Proof. } We omit the proof, which is similar to that of Theorem \ref{thm:newsob}.
$\hfill //$

\subsection{Refined heat kernel estimates}

In Theorems \ref{thm:mainthm} and \ref{thm:elec} we obtained heat kernel estimates
for operators $-\Delta -V$ where $V(x)=((N-2)/2)^2|x|^{-2}$ near the origin.
We shall now prove estimates for $-\Delta -V-V_1$, with $V_1$ also critical near
the origin. In Theorems \ref{thm:lognon} and \ref{thm:hk100} we consider the cases
$\Omega$ bounded and $\Omega=\R^N$ respectively.

For $k\geq 1$ and $\mu\leq 1/4$ we define
\begin{equation}
V_k^{\mu}(x)=\frac{\Bigl(\frac{N-2}{2}\Bigr)^2}{|x|^2}+\frac{1}{4|x|^2}\sum_{i=1}^{k-1}
X_1^2\ldots X_i^2  +\frac{\mu}{|x|^2}X_1^2\ldots X_k^2,\qquad x\in\Omega,
\label{eq:vkm}
\end{equation}
$(X_i=X_i(|x|/D),\; D=\sup_{\Omega}|x|\, )$
and consider the operator $H=-\Delta -V_k^{\mu}$ subject to Dirichlet boundary conditions
on $\partial\Omega$. In \cite[Proposition 7.2]{FT} the Hardy-Sobolev inequality
\begin{equation}
\int_{\Omega}(|\nabla u|^2-V_k^{1/4}u^2)dx
\geq c\Bigl(\int_{\Omega}|u|^{2N/(N-2)}(X_1\ldots X_{k+1})^{\frac{2N-2}{N-2}}
dx\Bigr)^{(N-2)/N}, \quad u\in C^{\infty}_c(\Omega),
\label{eq:5.1}
\end{equation}
was obtained.
Let $\beta$ be the largest solution of $\beta(1-\beta)=\mu$
and define
\begin{equation}
\phi_{k,\beta}(x)=|x|^{-\frac{N-2}{2}}X_1^{-1/2}\ldots X_{k-1}^{-1/2}X_k^{-\beta}.
\label{eq:5.3}
\end{equation}
Using (\ref{eq:dr}) we verify that $\Delta\phi_{k,\beta}+V_k^{\mu}\phi_{k,\beta}=0$ and hence
the change of variables $u=\phi_{k,\beta} w$ yields
\begin{equation}
\int_{\Omega}(|\nabla u|^2 -V_k^{\mu}u^2)dx=\int_{\Omega}|\nabla w|^2\phi_{k,\beta}^2 dx
\label{eq:5a}
\end{equation}
for all $w\in C^{\infty}_c(\Omega)$.
We have
\begin{theorem}
Let $\Omega$ be bounded, $1 \leq k <N-2$, and $0 <\mu \leq 1/4$. 
 The heat kernel of $H=-\Delta-V_k^{\mu}$ satisfies the estimate
\begin{equation}
K(t,x,y)<ct^{-N/2}\phi_{k,\beta}(x)\phi_{k,\beta}(y),\quad t>0,\; x,y\in\Omega;
\label{eq:5.02}
\end{equation}
here $V_k^{\mu}$ is given by (\ref{eq:vkm}) and   $\phi_{k,\beta}(x)$ by  (\ref{eq:5.3}).
\label{thm:lognon}
\end{theorem}
{\em Proof. }For the proof we distinguish two cases.\nl
{\em (1) Case $\mu<1/4$.} For $w\in C^{\infty}_c(\Omega)$ we have
\begin{eqnarray*}
\int_{\Omega}|\nabla w|^2\phi_{k,\beta}^2 dx &=&\int_{\Omega}|\nabla u|^2dx -
\int_{\Omega}V_k^{\mu}u^2dx\\
\mbox{(by (\ref{eq:5.01}))}\quad &\geq& c\left(\int_{\Omega}(|\nabla u|^2-V_{k-1}^{1/4}
u^2)dx\right) \\
\mbox{(by (\ref{eq:5.1}))}\quad
&\geq& c\Bigl(\int_{\Omega}|u|^{2N/(N-2)}(X_1\ldots X_k)^{\frac{2N-2}{N-2}}
dx\Bigr)^{(N-2)/N}\\
&=& c\Bigl(\int_{\Omega}|w|^{2N/(N-2)}|x|^{-N}(X_1\ldots X_{k-1})
X_k^{\frac{2N-2-2N\beta}{N-2}}
dx\Bigr)^{(N-2)/N}\\
&\geq& c\Bigl(\int_{\Omega}|w|^{2N/(N-2)}\phi_{k,\beta}^2
dx\Bigr)^{(N-2)/N}.
\end{eqnarray*}
This implies that $K_{\phi_{k,\beta}}(t,x,y)<ct^{-N/2}$ and (\ref{eq:5.02}) follows. \nl
{\em (2) Case $\mu=1/4$.} By \cite[Lemma 7.5]{FT} the following Sobolev inequality holds
\begin{eqnarray*}
&&\hspace{-1cm}\int_{\Omega}|\nabla w|^2|x|^{2-N}X_1^{-1}\ldots X_k^{-1}dx\geq \\
&\geq&c\Bigl(\int_{\Omega}|w|^{\frac{2N}{N-2}}|x|^{-N}X_1\ldots X_k 
X_{k+1}^{\frac{2N-2}{N-2}}dx\Bigr)^{(N-2)/N},
\quad w\in C^{\infty}_c(\Omega).
\end{eqnarray*}
This implies in particular
\[\int_{\Omega}|\nabla w|^2 \phi_{k,\beta}^2 dx\geq
c\Bigl(\int_{\Omega}|w|^{\frac{2N}{N-2}}\phi_{k,\beta}^2 dx\Bigr)^{(N-2)/N},
\quad w\in C^{\infty}_c(\Omega)\]
and hence we have the uniform estimate
$K_{\phi_{k,\beta}}(t,x,y)<ct^{-N/2}$ as required. //

We finally have the following consequence of Theorem \ref{thm:gian1},
where we retain the notation of that theorem: 
\begin{theorem}
Let   $1 \leq k <N-2$, and  $\epsilon\in(0,\epsilon_{k,0})$ with $\epsilon_{k,0}$ given by 
(\ref{eq:inf1}).   the heat kernel of
the operator $-\Delta-V_{k,\epsilon}$ satisfies
\[ K(t,x,y)<c t^{-N/2}\psi(x)\psi(y),\quad t>0,\; x,y\in\R^N;\]
\label{thm:hk100}
here $V_{k,\epsilon}$ is given by (\ref{eq:vke}) and  $\psi(x)$ by  (\ref{eq:defpsi}).
\label{thm:hk100}
\end{theorem}
{\em Proof.} The result follows directly from Theorem \ref{thm:gian1}
by means of (\ref{eq:eq}).$\hfill //$





\begin{thebibliography}{RRR}

\bibitem[AE]{AE}{Adimurthi, Esteban M. J. An Improved Hardy-Sobolev Inequality in $W^{1,p}$ and its
Application to Schr\"{o}dinger Operator. Preprint}

\bibitem[BV]{BV}{Brezis H. and V\'{a}zquez J.-L. Blow-up solutions
 of some nonlinear elliptic problems. Rev. Mat. Univ. Comp. Madrid
 {\bf 10} (1997) 443-469.}

\bibitem[CKN]{CKN}{Caffarelli L.A., Kohn R. V.  and Nirenberg L. First order
interpolation inequalities with weights. Compositio Math. {\bf 53} (1984)
259-275.}


 \bibitem[CM]{CM}{Cabr\'{e} X. and Martel Y. Existence versus instantaneous blowup for linear heat
 equations with singular potentials. C.R. Acad. Sci. Paris Ser. I Math. {\bf 329} (1999) 973-978.}


\bibitem[D]{D}{Davies E.B. Heat kernels and spectral theory. Cambridge University Press 1989.}

\bibitem[DD]{DD}{D\'{a}villa J. and Dupaigne L. Comparison Results for PDE's with a Singular
Potential. Preprint 2001.}

\bibitem[DS]{DS}{Davies E.B. and Simon B. $L^p$ norms of non-critical Schr\"{o}dinger semigroups.
J. Funct. Anal. {\bf 102} (1991) 95-115.}

\bibitem[FT]{FT}{Filippas S. and Tertikas A. Optimizing Improved
Hardy inequalities. J. Funct. Anal.  {\bf 192}  (2002)  186-233.}

\bibitem[G]{G}{Grigory'an A. Gaussian upper bounds for the
heat kernel on arbitrary manifolds. J. Diff. Geom. {\bf 45}
(1997) 33-52.}

\bibitem[LN]{LN}{Li Y. and Ni W.-M. On conformal scalar curvature equations in $\R^N$.
Duke Math. J. {\bf 57} (1988) 895-924.}

\bibitem[M]{M}{Maz'ya V. Sobolev spaces. Springer 1985.}

\bibitem[MS]{MS}{Milman P. and Semenov Yu.  Heat kernel bounds and 
 desingularizing weights. J.  Funct. Anal., to appear. }

\bibitem[T]{T}{Tertikas A. Critical phenomena in linear elliptic problems.
J. Funct. Anal. {\bf 154} (1998) 42-66.}

\bibitem[VZ]{VZ}{V\'{a}zquez  J.L. and Zuazua E. The Hardy inequality and the
asymptotic behaviour of the heat equation with an inverse-square
potential. J. Funct. Anal. {\bf 173} (2000) 103-153.}

\bibitem[WZ]{WZ}{Wang Zhi-Qiang  and Zhu Meijun Hardy Inequalities with Boundary Terms.
Preprint.}

\bibitem[Z]{Z}{Zhang Q.S. Global bounds of Schr\"{o}dinger heat kernels with negative potentials.
J. Funct. Anal. {\bf 182} (2001) 344-370.}

\end{thebibliography}
\end{document}